\documentclass[12pt]{amsart}
\usepackage{amssymb,amscd}
\usepackage{verbatim}
\usepackage{graphicx}

\newtheorem*{Whitney towers}{Theorem~\ref{Whitney towers}}
\newtheorem*{h-towers}{Theorems ~\ref{half} \& \ref{$(n)$-solvable}}

\newtheorem*{surgery curves}{Theorem~\ref{surgery curves}}
\newtheorem*{cg=0}{Theorem~\ref{vanish}}

\newtheorem{thm}{Theorem}[section]

\newtheorem{prop}[thm]{Proposition}

\theoremstyle{definition}

\newtheorem{note}[thm]{Note}

\numberwithin{equation}{section}
\numberwithin{figure}{section}

\newcommand{\x}{\times}
\newcommand{\np}{\newpage}
\newcommand{\nl}{\newline}
\newcommand{\Z}{\mathbb{Z}}

\newcommand{\Q}{\mathbb{Q}}

\newcommand{\R}{\mathbb{R}}

\newcommand{\f}{\noindent}

\begin{document}
\title{Local move formulae for \\the Alexander polynomials of  $n$-knots }
\author{Eiji Ogasa\\
}

\begin{abstract} 
It is well-known:Suppose there are 
three 1-dimensional links $K_+$, $K_-$, $K_0$ 
such that  $K_+$, $K_-$, and $K_0$ coincide out of 
a 3-ball $B$ trivially embedded in $S^3$ 
and that   $K_+\cap B$, $K_-\cap B$, and $K_0\cap B$ are drawn as follows. 
Then $\Delta_{K_+}-\Delta_{K_+}=(t-1)\cdot\Delta_{K_0}$, 
where $\Delta_{K}$ is the Alexander polynomial of $K$. 

$${\mathrm{Figure 1}}$$  

\noindent 
We know similar formulae of other invariants of 1-dimensional knots and links.
(The Jones polynomial etc.)

It is natural to ask:   
Suppose there are 
two $n$-dimensional knots $K_+$, $K_-$ and a submanifold $K_0$ 
such that  $K_+$, $K_-$, and $K_0$ coincide out of 
a $n$-ball $B$ trivially embedded in $S^{n+2}$.  
Then is there a relation in $K_+\cap B$, $K_-\cap B$, and $K_0\cap B$ with the following property(*)? 
(*)If $K_+$, $K_-$, and $K_0$ satisfy this relation, 
an invariant of $K_+$, that of $K_-$, and that of $K_0$ 
satisfy a fixed relation. 

In this paper we pove there are such a relation where  
$K_+$, $K_-$, and $K_0$ satisfy the formula 
$\Delta_{K_+}-\Delta_{K_+}=(t-1)\cdot\Delta_{K_0}$, 
where $\Delta_{K}$ is a polynomial to represent the Alexander polynomial of $K$.

We show another relation where 
$K_+$, $K_-$, and $K_0$ satisfy the formula 

\f
${\mathrm{Arf}}K_+-{\mathrm{Arf}}K_-=\{|bP_{4k+2}\cap I(K_0)|+1\}mod 2,$

\f
where (1)$I(\quad)$ is the inertia group. and $I(K_0)$ is the inertia group of a smooth manifold which is orientation preserving diffeomorphic to $K_0$. 
(2)For a group $G$, $|G|$ denote the order of $G$. 

A local move formula is a relation of an invariant of a few knots 
related by a local move as above.

\end{abstract}

\thanks{This research was partially supported by Research Fellowships
of the Promotion of Science for Young Scientists.
\nl{\bf PACS nos.} 11-25w, 11-25Uv.}
\maketitle

\section{Introduction} \label{intro} 

It is well-known:Suppose there are 
three 1-dimensional links $K_+$, $K_-$, $K_0$ 
such that  $K_+$, $K_-$, and $K_0$ coincide out of 
a 3-ball $B$ trivially embedded in $S^3$ 
and that   $K_+\cap B$, $K_-\cap B$, and $K_0\cap B$ are drawn as follows. 
Then $\Delta_{K_+}-\Delta_{K_+}=(t-1)\cdot\Delta_{K_0}$, 
where $\Delta_{K}$ is the Alexander polynomial of $K$.

$${\mathrm{Figure 1}}$$  

\noindent 
We know similar formulae of other invariants of 1-dimensional knots and links.
(The Jones polynomial etc.
See \S5 of \cite{Kauffman}. See also \cite{Cn} \cite{J}.)

It is natural to ask:   
Suppose there are 
two $n$-dimensional knots $K_+$, $K_-$ and a submanifold $K_0$ 
such that  $K_+$, $K_-$, and $K_0$ coincide out of 
a $n$-ball $B$ trivially embedded in $S^{n+2}$.  
Then is there a relation in $K_+\cap B$, $K_-\cap B$, and $K_0\cap B$ with the following property(*)? 
(*)If $K_+$, $K_-$, and $K_0$ satisfy this relation, 
an invariant of $K_+$, that of $K_-$, and that of $K_0$ 
satisfy a fixed relation. 

In this paper we pove there are such a relation where  
$K_+$, $K_-$, and $K_0$ satisfy the formula 
$\Delta_{K_+}-\Delta_{K_+}=(t-1)\cdot\Delta_{K_0}$, 
where $\Delta_{K}$ is a polynomial to represent the Alexander polynomial of $K$.

We show another relation where 
$K_+$, $K_-$, and $K_0$ satisfy the formula 

\f
${\mathrm{Arf}}K_+-{\mathrm{Arf}}K_-=\{|bP_{4k+2}\cap I(K_0)|+1\}mod 2,$

\f
where (1)$I(\quad)$ is the inertia group. and $I(K_0)$ is the inertia group of a smooth manifold which is orientation preserving diffeomorphic to $K_0$. 
(2)For a group $G$, $|G|$ denote the order of $G$. 

A {\it{local move formula}} is a relation of an invariant of a few knots related by a local move as above.

\cite{Ogasa1} is a preprint of this paper. 

The author proved another local move formulae in \cite{O1}\cite{O2}.

\section{Review of the Alexander polynomials for $n$-knots and $n$-links} 
\label{alex}  

We review 
the Alexander polynomials for $n$-knots and $n$-links and $n$-submanifolds. 
See \cite{Alexander} \cite{L1}, \cite{Ml}, \cite{R} for detail.

\noindent
We work in the smooth category. 
Let $K=(K_1,...,K_m)$ be an $n$-dimensional closed oriented submanifold of $S^{n+2}$. 
It is known any tubular neighborhood of $K$ is 
$K\x D^2$. (See P.49, 50 of \cite{Kr}.)
Put $X=\overline{S^{n+2}-K\x D^2}$. 
Then any $S^1$ in $X$ is oriented 
by using the orientation of $S^{n+2}$ and 
that of $K$.
Let $\iota:S^1\rightarrow X$ denote the embedding. 
Take a homomorphism 
$\alpha:H_1(X;\Z)\rightarrow\Z$ 
such that  
$\alpha\circ\iota_*:H_1(S^1;\Z)
\rightarrow\Z$ carries +1 to +1. 
Then the infinite cyclic covering 
$\pi:$ $\widetilde X\rightarrow X$ 
associated with $\alpha$ is called 
the {\it cannonical cyclic covering} of $K$. 

\noindent 
We can regard $H_p(\widetilde X;\Z)$ as a $\Z[t,t^{-1}]$-module 
by using 
the covering translation 
$\widetilde X\rightarrow\widetilde X$ 
defined by $\alpha$. 

\noindent 
We can also regard $H_p(\widetilde X;\Q)$ as a $\Q[t,t^{-1}]$-module.

\noindent 
Module theory says that any $\Q[t,t^{-1}]$-module  is congruent to 

\f $(\Q[t,t^{-1}]/{\lambda_1})\oplus...\oplus(\Q[t,t^{-1}]/{\lambda_l})
\oplus(\oplus^k\Q[t,t^{-1}])$, 
where $\lambda_*\in$ $\Q[t,t^{-1}]$ is not zero and 
$\lambda_*$ is not the $\Q[t,t^{-1}]$-balanced class of 1. 

\f
Two polynomials $f(t), g(t)\in\Q[t,t^{-1}]$ are said to be 
{\it $\Q[t,t^{-1}]$-balanced } 
(written $f\doteq g$) if there is an integer $n$ 
and a nonzero rational number $r$ 
such that 
$f(t)=r\cdot t^n\cdot g(t)$.

\noindent
Let $H_p(\widetilde X;\Q)$ be 
$\Bbb Q[t,t^{-1}]/{\lambda_1}\oplus...
\oplus\Bbb Q[t,t^{-1}]/{\lambda_l}\oplus^k\Bbb Q[t,t^{-1}]$ as above. 
{\it The $\Bbb Q[t,t^{-1}]$-$p$-Alexander polynomial} is 
the $\Bbb Q[t,t^{-1}]$-balanced class of 
the product $\lambda_1\cdot...\cdot\lambda_l$  
if $k=0$, 
where $k$ is the rank of the free part.     
{\it The $\Bbb Q[t,t^{-1}]$-$p$-Alexander polynomial} is 0 if $k\neq0$. 
If $H_p(\widetilde X;\Q)\cong0$,      
{\it The $\Bbb Q[t,t^{-1}]$-$p$-Alexander polynomial} is 1. 

We discuss the  $\Q[t,t^{-1}]$ module case 
but our results can be extended to the  $\Z[t,t^{-1}]$ module case. 

In this paper we mainly discuss the case where $K$ is a knot although 
we discuss other cases a little. 
Furthermore, our results can be extended to 
some other cases without heavy difficulty.

If $K$ above is a connected smooth manifold 
which is PL homeomorphic to the standard sphere, 
$K$ is called {\it $n$-(dimensional )knot}   
(See \cite{CochranOrr} etc).

\section{Main results}\label{main}  

In this section, we prove local move formulae 
for $n$-knots  $\subset S^{n+2}$. ($n\geqq3$.)

Let $K_+$, $K_-$ be an $n$-knot $\subset S^{n+2}$ ($n\geqq3$). 
Let $K_0$ be an $n$ submanifold $\subset S^{n+2}$. 
Let $B$ be an $(n+2)$-ball trivially embedded in $S^{n+2}$. 
Suppose that $K_+$ coincides with  $K_-$ in $S^{n+2}-B$. 
Note that there is a Seifert hypersurface 
$V_+$ (resp. $V_-$) for 
$K_+$ (resp. $K_-$) 
such that 
$V_+$ coincides with  $V_-$ in $S^{n+2}-B$. 
Suppose that 
$V_+\cap B$  (resp. $V_-\cap B$) is a disjoint union of 
an $(n+1)$-dimensional $p$-handle 
$h^{p}_+$  (resp. $h^{p}_-$) 
and an $(n+1)$-dimensional $(n+1-p)$-handle 
$h^{n+1-p}_+$  (resp. $h^{n+1-p}_-$) 
  which are attached to $\partial B$ and 
  which are embedded trivially in $B$. 
Let $p\neq n+1-p$. 
Suppose that 
$h^{n+1-p}_+$ coincides with $h^{n+1-p}_-$. 
Suppose that the linking number (in $B$) of 
`$h^{p}_+\cup (-h^{p}_-)$' 
and 
`$h^{n+1-p}_+$ whose attached part is fixed in $\partial B$' is one 
if an orientation is given.  The concept is drawn in Figure 2. 
Let $K_0$ be $\partial(V_+ -{\mathrm {Int}}B)$. 
\hskip3cm [Figure 2]

\noindent
Then we say that  $(K_+$, $K_-$, $K_0)$ is related 
by the {\it $(p,n+1-p)$-move.}

We draw the figure of the $(1,2)$-move case (the case if $p=1$ and $n=2$) 
in Note below Theorem \ref{surface}.

Let $n=4k+1$ in the above case. 
Suppose that 
$K_+$, $K_-$, $V_+$, $V_-$ satisfy the same condition at $S^{4k+3}-B$ 
as in (i). 
Suppose $V_+\cap B$ (resp. $V_-\cap B$) is 
a $(4k+2)$-dimensional $(2k+1)$-handle $h_+$ (resp. $h_-$). 
Supoose that the core of $h_+$ (resp. $h_-$) is trivially embedded in $B$. 
Push off the core in the positive direction of the normal bundle of 
$V_+$ (resp. $V_-$) in $B$. 
Note that we can consider the framing (in $B$) of $h_+$ (resp. $h_-$) . 
Suppose that the framing of $h_+$ (resp. $h_-$) is 0 (resp. 1)
if an orientation is given. 
Let $K_0$ be $\partial(V_+ -{\mathrm {Int}}B)$. 
(The 1-dimensional case of 
this relation among $K_+, K_-, K_0$ is one in Figure 1.)

\f 
Then we say that  $(K_+$, $K_-$, $K_0)$ is related 
by the {\it $XXII$-move.}

\begin{note}\label{meaning}
{\rm{
One way of saying, when we make $K_-$, $K_0$, 
we just operate in $B$ and we do not need 
the diffeomorphism type or the hmeomorphism type 
of $K_-$, $K_0$. 
In this meaning, we use the word `local' in the above definition. 
}}
\end{note}

\begin{thm}\label{knots}
Let $K_+$, $K_-$ be $n$-knots $\subset S^{n+2}$ ($n\geqq3$). 
Let $K_0$ be $n$-submanifold $\subset S^{n+2}$. 
Suppose that  $(K_+$, $K_-$, $K_0)$ is related by the $(p,n+1-p)$-move.  
Then we have: 

\f
$\Delta^p_{K_+}-\Delta^p_{K_+}=(t-1)\cdot\Delta^p_{K_0}$, 

\f
where  
$\Delta^p_{K}$ is a polynomial 
whose balanced class is the $p$-Alexander polynomial for $K$.



\end{thm}


\begin{thm}\label{middle}
Let $K_+$, $K_-$ be $(4k+1)$-knots. 
Let $K_0$ be a closed oriented $(4k+1)$-submanifold $\subset S^{4k+3}$.
Suppose that  $(K_+$, $K_-$, $K_0)$ is related by the XXII-move.  
Then we have: 

\f
$\Delta^{(2k+1)}_{K_+}(t)-\Delta^{(2k+1)}_{K_-}(t)=
(t-1)\cdot\Delta^{(2k+1)}_{K_0}(t)$.  

\f
where  
$\Delta^{(2k+1)}_{K}$ is a polynomial 
whose balanced class is the $p$-Alexander polynomial for $K$.
\end{thm}

\f{\bf{Note.}} 
For $(4k+3)$-knots, we can define XXII-move.  
However, note the following: 
Suppose $K_-$and $K_0$ satisfy the XXIIrelation for a  $(4k+3)$-knot $K_+$. 
Then $K_0$ is not a knot in general. 
Because: there is an example such that 
$K_+\cong S^3$ and $K_-\cong\R P^3$.


\begin{thm}\label{inertia}
Let $K_+$, $K_-$, $K_0$, $B$ be as in Theorem \ref{middle}. 
Let $bP_{4k+2}$ be the $bP$-subgroup $\subset \Theta^{4k+1}$.
Suppose $bP_{4k+2}$ is not congruent to the trivial group. 
Then we have 

\f
${\mathrm{Arf}}K_+-{\mathrm{Arf}}K_-=\{|bP_{4k+2}\cap I(K_0)|+1\}mod 2,$

\f
where 
(1)$I(\quad)$ is the inertia group. $I(K_0)$ is the inertia group of a smooth manifold which is orientation preserving 
diffeomorphic to $K_0$. 

\f
(2)For a group $G$, $|G|$ denote the order of $G$. 
\end{thm}

\f{\bf{Note.}}
See [KM] for the $bP$-subgroup. 
See [Kk] [BS] for the inertia group.


\f{\bf Proof of Theorem \ref{knots} and \ref{middle}. }
Let $V_*$ be a compact oriented $(n+1)$-submanifold $\subset S^{n+2}$ 
such that  $\partial V_*=K_*$ (their orientaion are compatible). 
Recall $V_*$ is called a Seifert hypersurface for $K_*$.
\f
In Theorem \ref{middle} we put $n=4k+1$ and $p=2k+1$. 

\f
Take $X_*$ $\widetilde X_*$ as in \S\ref{alex}. 
Let $Y_*=X-V_*\x[-1,1]$, $V_*\x[-1,1]$ is the tubular neibourhood of 
$V_*$ in $Y$. 
Consider the Meyer-Vietoris exact sequence: 

\f
$
H_*(\amalg_{-\infty}^{\infty}V_*)\stackrel{f_*}\to
H_*(\amalg_{-\infty}^{\infty}Y_*)\to
H_p(\widetilde X_*;\Z)
$.  

\f 
There are $V_*$ such that $f_*$ is represented by the following matrixes:
$$
P_+=\{p^+_{ij}\}, P_-=\{p^-_{ij}\}, P_0=\{p^0_{ij}\}, 
$$
such that 
(1) $P_+$ and $P_-$ are $(n+1)\x(n+1)$ matrices. 
$P_0$ is an $n\x n$ matrix.
(2)$p^+_{n+1,n+1}-p^-_{n+1,n+1}=t-1,$
(3)$p^+_{ij}=p^-_{ij}=p^0_{ij}$. ($1\leqq i\leqq n$, $1\leqq j\leqq n, $). 

\f
By calculus of determinants, 
${\mathrm{det}}P_+-{\mathrm{det}}P_-=(t-1){\mathrm{det}}P_0.$

\f 
Module theory says that 
${\mathrm{det}}P_*$ represents the $p$-Alexander polynomial for $K_*$. 
Hence 
Theorem \ref{knots} and \ref{middle} hold.

\f{\bf Proof of Theorem \ref{inertia}.}
By [KM], $bP_{4k+2}=\Bbb Z_2$. 
In our case, the Arf invariant of a knot coincides with 
that of a manifold diffeomorphic to the knot. 
(See [Br1] for more results  on $bP_{4k+2}$. )
Put $bP_{4k+2}=\{1,g\}$, where $g^2=1.$

\f 
Let $V$ be the total space of $D^{2k+1}$ bundle over $S^{2k+1}$ 
associated with the tangent bundle of $S^{2k+1}$. 
Let $\bigvee_p$ denote a plumbing ( see [Br2] ). 
Then $\partial(V\bigvee_p V)$ represents $g\in bP_{4k+2}$.

\f
Put $M=\partial V$.
By \cite{JW} and \cite{Adams}, 
$M$ is homotopy type equivalent to $S^{2k+1}\x S^{2k}$
if and only if $k=0,1,3.$ 
Hence $M$ is not diffeomorphic to $S^{2k+1}\x S^{2k}$ in our case.

\f
\cite{BS} proved $M\sharp\Sigma=M$. ( 
Hence 
$I(M)\cap bP_{4k+2}=$ 
$I(M\sharp\Sigma)\cap bP_{4k+2}=\Bbb Z_2$. 

\f 
Corollary 3 of [Kk] proved $I(S^p\x S^q)=\{1\}$ for $p+q\geqq5$. 


\f
By [Kk], 
$S^{2k+1}\x S^{2k}\sharp \Sigma$ is not diffeomorphic to $S^{2k+1}\x S^{2k}$. 
Hence $I(S^{2k+1}\x S^{2k}\sharp \Sigma)\cap bP_{4k+2}$=$\{1\}$.

There are four cases. 
Put $T=S^{2k+1}\x D^{2k+1}$. 
Let $\cong$ denote a diffeomorphism. 

\f  
(1)$K_+\cong\partial(V\bigvee_p V)$,  
   $K_-\cong\partial(V\bigvee_p T)$,  
     $K_0\cong\partial V$.

\f 
(2)$K_+\cong\partial(T\bigvee_p V)$, 
   $K_-\cong\partial(T\bigvee_p T)$ 
    $K_0\cong\partial T$.  

\f
(3)$K_+\cong\partial\{(V\bigvee_p V) \natural (V\bigvee_p V)\}$, 
   $K_-\cong\partial\{(V\bigvee_p T)\natural (V\bigvee_p V)\}$,  
    $K_0\cong\partial \{V\natural (V\bigvee_p V)\}$.  

\f
(4)$K_+\cong\partial\{(T\bigvee_p V) \natural (V\bigvee_p V)\}$, 
   $K_-\cong\partial\{(T\bigvee_p T) \natural (V\bigvee_p V)\}$,  
   $K_0\cong\partial\{T\natural (V\bigvee_p V))$.

\f 
The formula in Theorem \ref{inertia} 
holds in each case by the above discussions. 
Hence the formula holds.

\section{More results in the 2-knot case}\label{two}
Our main results can be extended to some other cases where 
$K_+$ (resp. $K_-$) is not a knot. 
In this section we show more results 
in the case of 2-dimensional sbmanifold case.

\begin{thm}\label{surface}    
Let $\Sigma_1,...,\Sigma_{\alpha}$ be  connected closed oriented surfaces. 
Let $g_i$ be the genus of $\Sigma_i$. 
Put $\beta=\Sigma_{i=1}^{\alpha}g_i$. 
Let $K_+$ (resp. $K_-$) be a 2-dimensional submanifold $\subset S^4$ 
which is diffeomorphic to a disjoint ordered oriented manifold 
$(\Sigma_1,...,\Sigma_{\alpha})$. 
Suppose $\alpha=\beta+1$.
Let $K_0$ be a 2-dimensional submanifold $\subset S^4$.  
Suppose that  $(K_+, K_-, K_0)$ is related by the $(1,2)$-move. 
Then we have:
$$\Delta_{K_+}(t)-\Delta_{K_-}(t)=(t-1)\cdot\Delta_{K_0}(t).$$  
\end{thm}

\noindent{\bf Note.} (1)
If $\alpha=1$, $K_+$ and $K_-$ are $S^2$-knots. 
Then $K_0$ is diffeomorphic to $S^2$ or $S^2\amalg T^2$.
In each case Theorem 1 holds. 
In general, if we put 
$H_0(K_0;\Bbb Q)\cong\Bbb Q^{\alpha'}$ 
and 
$H_1(K_0;\Bbb Q)\cong\Bbb Q^{2\beta'}$, 
then $\alpha'=\beta'+1$.

\f
(2) 
Since $(K_+, K_-, K_0)$ is related by the $(1,2)$-move, 
there is a 4-ball $B$ trivially embedded in $S^4$ 
with the following properties. 
We regard $B$ as (2-disc)$\x[0,1]\x\{t\vert -1\leqq t\leqq1\}$.  

\f
(i) $K_+-B$, $K_--B$, and $K_--B$ coincide each other.   

\f
(ii) $B\cap K_+$,  $B\cap K_-$, $B\cap K_-$ 
are drawn as in Figure 3.

\f 
In Figure 3 we draw 
$B_{-0.5}\cap K_*$, $B_0\cap K_*$, $B_{0.5}\cap K_*$, 
where $B_{t_0}$=(2-disc)$\x[0,1]\x\{t\vert t=t_0 \}$.  
We suppose that 
each vector $\overrightarrow{x}$, $\overrightarrow{y}$ 
in Figure 3 is a tangent vector of each disc at a point. 
(Note we use $\overrightarrow{x}$ (resp. $\overrightarrow{y}$)
for different vectors.)
The orientation of each disc 
in Figure 3 is determined by the each set 
$\{\overrightarrow{x},\overrightarrow{y}\}$.

\f 
In \cite{O1} 
the author calls the operation to change $K_+$ into $K_-$ (1,2)-pass-move. 
Around Figure 4.1 and 4.2 in \cite{O1}, 
we wrote more explanation of 
the figure of $B\cap K_+$ and that of $B\cap K_-$.

\f
(4) 
After sending these results (without Appendix) to several people, 
the author was informed Giller's article, P.627,628 of [Gi].  
Only in the $n=2$ case  
Giller proved a result which is weaker  than ours. 
See the Appendix for detail. 

\hskip5cm[Figure 3]%

\vskip3mm
\f{\bf{Proof of Theorem \ref{surface}.}} 
Since $\beta=\Sigma_{i=1}^{\alpha}g_i$, 
all Seifert surfaces $V_*$ for $K_*$ 
has a property that 
$H_1(V_*;\Q)\cong H_2(V_*;\Q)$ 
and that 
$H_0(V_*;\Q)\cong\Q.$

\f 
The left of the proof is same as the proof of Theorem \ref{knots}. 
\vskip3mm

On the condition $\alpha=\beta+1$ in Theorem \ref{surface} we have: 

\begin{prop}\label{rest}
We CANNOT remove the condition $\alpha=\beta+1$ 
in Theorem \ref{surface} in general.  
\end{prop}

\f{\bf{Proof of Theorem \ref{rest}.}}
Let $K_-$, $K_+$, $K_0$ be 2-dimensional oriented subamanifold
$\subset S^4$ 
which are diffeomorphic to $T^2$-knots.  
Suppose $(K_-, K_+, K_0)$ is related by the $(1,2)$-move.

Suppose that 

\f(1) $K_+$ bounds 
 $V_+\cong S^1\x B^2\natural (S^2\x S^1-B^3)\natural (S^2\x S^1-B^3)$

\f(2) $K_-$ bounds $V_-\cong V_+$. 

\f(3) $K_0$ bounds $V_0\cong S^1\x B^2\natural (S^2\x S^1-B^3)$. 

Consider the exact sequence as in Proof of Theorem \ref{knots}:

\f
$
H_*(\amalg_{-\infty}^{\infty}V_*)\stackrel{f_*}\to
H_*(\amalg_{-\infty}^{\infty}Y_*)\to
H_p(\widetilde X_*;\Z)
$.

We can suppose that 

\f(1) 
$f_+$ is represented by 
$\begin{pmatrix}
3t-2&0&0\\
t-1&2t-1&0\\
\end{pmatrix}$.

\f(2) 
$f_-$ is represented by 
$\begin{pmatrix}
3t-2&0&t-1\\
t-1&2t-1&0\\
\end{pmatrix}$.

\f(3) 
$f_0$ is represented by 
$\begin{pmatrix}  
t-1&2t-1\\
\end{pmatrix}$.

The above exact sequences are:
$$\CD
\oplus^3\Bbb Z[t,t^{-1}]@>>
{\begin{pmatrix}
3t-2&0&0\\
t-1&2t-1&0\\
\end{pmatrix}}
>\oplus^2\Bbb Z[t,t^{-1}]@>>>
\Bbb Z[t,t^{-1}]/\{(3t-2)\}\oplus\Bbb Z[t,t^{-1}]/\{(2t-1)\}@>>>0
\endCD$$

$$\CD
\oplus^3\Bbb Z[t,t^{-1}]@>>
{\begin{pmatrix}
3t-2&0&t-1\\
t-1&2t-1&0\\
\end{pmatrix}}
>\oplus^2\Bbb Z[t,t^{-1}]@>>>0@>>>0
\endCD$$

$$\CD
\oplus^2\Bbb Z[t,t^{-1}]@>>
{\begin{pmatrix}
t-1&2t-1\\
\end{pmatrix}}
>\Bbb Z[t,t^{-1}]@>>>0@>>>0
\endCD$$ 


\f
Therefore we have: 

\f (1) $(3t-2)(t-1)$ represents the Alexander polynomial of $K_+$. 

\f(2) $1$ represents the Alexander polynomial of $K_-$. 

\f(3) $1$ represents the Alexander polynomial of $K_0$.

For the above $K_-$, $K_+$, $K_0$, 
$\Delta_{K_+}-\Delta_{K_-}=(t-1)\cdot\Delta_{K_0}$ 
DOES NOT hold 
for any set of polynomials 
$\Delta_{K_-}$, $\Delta_{K_+}$, $\Delta_{K_0}$ 
such that 
$\Delta_{K_-}$ (resp. $\Delta_{K_+}$, $\Delta_{K_0}$) 
represents the  Alexander polynomials of 
$K_-$ (resp. $K_+$, $K_0$). 
The proof is completed.

\vskip3mm
Next we discuss ``normalization'' of the Alexander polynomials. 
Recall that, in the case of 1-links, 
we can choose a unique polynomial from the all polynomilas 
whose ballanced classes are the Alexander polynomial.  
(See e.g. \cite{Kauffman} for detail.) 
However, we have:

\begin{prop}\label{nonnor}
In the $n\neq1$ case in Theorem \ref{knots},
we CANNOT choose a unique polynomial 
from all polynomials which represent the Alexander polynomial 
to be compatible with our local move formula.
\end{prop}

We can suppose that $K_+$, $K_-$, $K_0$ are trivial knots 
and that $(K_+, K_-, K_0)$ can be related by the $(1,2)$-move. 
Becase: 
Let $V_+\cong\overline{S^2\x S^1-B^3}$. 
Let $V_-\cong\overline{S^2\x S^1-B^3}$. 
Let $V_0\cong\overline{B^3}$. 
Use these $V_+$, $V_-$, $V_0$. 

If we can take a unique polynomial to represent the Alexander polynomial, 
then we can let the Alexander polynomial 
$a\cdot t^m$ for $K_+$, $K_-$ and $K_0$, 
where $a$ is a nonzero rational number.   
Hence $a\cdot t^m-a\cdot t^m=(t-1)\cdot a\cdot t^m$. 
Hence $0=(t-1)\cdot a\cdot t^m$. 
It is the contardiction.  
Hence we CANNOT choose a unique polynomial.

\vskip3mm
\f 
Computer Science, Meijigakuin University, Yokohama, Kanagawa, 244-8539, Japan 

\f pqr100pqr100@yahoo.co.jp

\vskip5mm
\hskip3cm\f{\bf{Appendix.}}

\f 
We explain the fact in Theorem \ref{surface}(4) a little more.

\f
Replace Figure 3 with Figure 3 
in the definition that 
$(K_+, K_-, K_0)$ is related by the $(1,2)$-move. 
Then we say that $(K_+, K_-, K_0)$ is {\it related by the ribbon-move.} 

\f{\bf Note.}  
In \cite{O1} the author call the operation to change $K_+$ into $K_-$ 
(resp. $K_-$ into $K_+$) in Figure 4, 
ribbon-move.

\f
By using Theorem \ref{surface} in this paper and 
Proposition 4.2, 4.3, 4.4 in \cite{O1},
we have: 
If $(K_+, K_-, K_0)$ is  related by the ribbon-move,  
then 
$(K_+, K_-, K_0)$ is  related by the (1,2)-move.

\f Thus  we have: 

\f{\bf{Theorem.}}
{\it 
In Theorem \ref{surface}, 
we can replace the word  `the (1,2)-move' with `the ribbon-move.' 
}

\f{\bf Note.} 
In P.627, 628 of [Gi], 
Giller proved a weaker case of this Theorem: 
[Gi] does not prove 
the case where $K_+$ is a sphere, $K_-$ is a sphere, 
and $K_0$ is not a sphere. 
It means that [Gi]'s formula is not a local move formula 
in the meaning of Note \ref{meaning}. 
Furthermore [Gi] does not prove 
more than one of $K_+$, $K_-$, $K_0$ is a sphere. 
In the meaning of the following Proposition, 
ours are stronger than the formula in [Gi].

\hskip5cm[Figure 4]

\f{\bf Note.}
The orientation of the part of $B\cap K_0$ derived from $B\cap K_+$ (resp.  $B\cap K_+$) 
is given by using $B\cap K_+$ (resp.  $B\cap K_+$).  
The orientation of $B\cap K_0$ is compatible with 
the part of $B\cap K_0$ derived from $B\cap K_+$ (resp.  $B\cap K_+$).

It is natural to ask the following. 
If $(K_+, K_-, K_0)$ is related by the $(1,2)$-move, 
then do they compose a triple of Figure 1.1? 

\f The answer is negative in general by the following Proposition.

Let $K=(K_1, K_2, K_3)$ and $K'=(K'_1, K'_2, K'_3)$ 
be 2-dimensional submanifolds $\subset S^4$ 
such that 
$K_1\cong K'_1\cong S^2$,  
$K_2\cong K'_2\cong S^2$, 
and  
$K_3\cong K'_3\cong \Sigma_2$, 
 $\Sigma_2$ is the oriented closed surface with the genus two.  

\f 
Supoose that 
alk$(K_1, K_3)$ is one, 
alk$(K_2, K_3)$ is one, 
where alk$(\quad)$ denotes the alinking number (in \cite{S}). 

\f 
Suppose that 
$K$ is changes into $K'$ by one (1,2)-pass-move (see \cite{O1}). 
We can suppose this (1,2)-pass-move 
let alk$( K_1, K_3 )$ zero 
and 
let alk$( K_2, K_3 )$ zero.

Then we prove: 

\f{\bf{Proposition.}}{\it 
Let $K$ and $K'$ be as above. 
$K$ does not change into $K'$ by one ribbon-move. 
}

\f{\bf{Proof.}}  
One ribbon-move cannot change  
alk$( K_1, K_3 )$ and 
alk$( K_2, K_3 )$ together.



\np 
\unitlength 0.1in
\begin{picture}(46.22,13.80)(8.10,-21.90)
%
\special{pn 8}%
\special{ar 1356 1356 546 546  0.8902751 6.2831853}%
\special{ar 1356 1356 546 546  0.0000000 0.8502422}%
%
\special{pn 8}%
\special{pa 946 1716}%
\special{pa 1736 986}%
\special{fp}%
\special{sh 1}%
\special{pa 1736 986}%
\special{pa 1673 1017}%
\special{pa 1697 1022}%
\special{pa 1701 1046}%
\special{pa 1736 986}%
\special{fp}%
%
\special{pn 8}%
\special{pa 1266 1276}%
\special{pa 976 976}%
\special{fp}%
\special{sh 1}%
\special{pa 976 976}%
\special{pa 1008 1038}%
\special{pa 1013 1014}%
\special{pa 1037 1010}%
\special{pa 976 976}%
\special{fp}%
%
\special{pn 8}%
\special{pa 1416 1416}%
\special{pa 1726 1746}%
\special{fp}%
%
\special{pn 8}%
\special{ar 3086 1396 546 546  0.8902751 6.2831853}%
\special{ar 3086 1396 546 546  0.0000000 0.8502422}%
%
\special{pn 8}%
\special{pa 2996 1316}%
\special{pa 2706 1016}%
\special{fp}%
\special{sh 1}%
\special{pa 2706 1016}%
\special{pa 2738 1078}%
\special{pa 2743 1054}%
\special{pa 2767 1050}%
\special{pa 2706 1016}%
\special{fp}%
%
\special{pn 8}%
\special{pa 3146 1456}%
\special{pa 3456 1786}%
\special{fp}%
%
\special{pn 8}%
\special{pa 3210 1520}%
\special{pa 2990 1320}%
\special{fp}%
\special{pa 3240 1550}%
\special{pa 3100 1420}%
\special{fp}%
%
\special{pn 8}%
\special{pa 3140 1330}%
\special{pa 3480 1010}%
\special{fp}%
\special{sh 1}%
\special{pa 3480 1010}%
\special{pa 3418 1041}%
\special{pa 3441 1047}%
\special{pa 3445 1070}%
\special{pa 3480 1010}%
\special{fp}%
%
\special{pn 8}%
\special{pa 3010 1470}%
\special{pa 2690 1770}%
\special{fp}%
%
\special{pn 8}%
\special{ar 4886 1396 546 546  0.8902751 6.2831853}%
\special{ar 4886 1396 546 546  0.0000000 0.8502422}%
%
\special{pn 8}%
\special{pa 4570 1830}%
\special{pa 4540 970}%
\special{fp}%
\special{sh 1}%
\special{pa 4540 970}%
\special{pa 4522 1037}%
\special{pa 4542 1023}%
\special{pa 4562 1036}%
\special{pa 4540 970}%
\special{fp}%
%
\special{pn 8}%
\special{pa 5240 1820}%
\special{pa 5230 970}%
\special{fp}%
\special{sh 1}%
\special{pa 5230 970}%
\special{pa 5211 1037}%
\special{pa 5231 1023}%
\special{pa 5251 1036}%
\special{pa 5230 970}%
\special{fp}%
\put(11.2000,-22.2000){\makebox(0,0)[lb]{K}}%
\put(30.5000,-22.7000){\makebox(0,0)[lb]{K}}%
\put(49.0000,-22.9000){\makebox(0,0)[lb]{K}}%
\put(12.7000,-23.2000){\makebox(0,0)[lb]{+}}%
\put(31.7000,-23.3000){\makebox(0,0)[lb]{-}}%
\put(50.0000,-23.6000){\makebox(0,0)[lb]{0}}%
\end{picture}%
\vskip2cm
Figure 1.

\np 
\unitlength 0.1in
\begin{picture}(53.62,22.22)(4.10,-31.30)
%
\special{pn 8}%
\special{ar 1820 1830 922 922  0.8926770 6.2831853}%
\special{ar 1820 1830 922 922  0.0000000 0.8621701}%
%
\special{pn 8}%
\special{pa 1100 1270}%
\special{pa 2260 2630}%
\special{fp}%
%
\special{pn 8}%
\special{pa 2550 2380}%
\special{pa 1400 1020}%
\special{fp}%
%
\special{pn 8}%
\special{pa 2260 1020}%
\special{pa 1850 1410}%
\special{fp}%
\special{pa 1470 1810}%
\special{pa 1000 2250}%
\special{fp}%
\special{pa 2540 1280}%
\special{pa 2090 1680}%
\special{fp}%
\special{pa 1690 2100}%
\special{pa 1250 2520}%
\special{fp}%
%
\special{pn 8}%
\special{ar 4850 1870 922 922  5.5896934 6.2831853}%
\special{ar 4850 1870 922 922  0.0000000 5.5599667}%
%
\special{pn 8}%
\special{pa 4301 2599}%
\special{pa 5643 1418}%
\special{fp}%
%
\special{pn 8}%
\special{pa 5389 1132}%
\special{pa 4047 2303}%
\special{fp}%
%
\special{pn 8}%
\special{pa 4033 1443}%
\special{pa 4430 1847}%
\special{fp}%
\special{pa 4835 2221}%
\special{pa 5282 2684}%
\special{fp}%
\special{pa 4289 1159}%
\special{pa 4696 1603}%
\special{fp}%
\special{pa 5122 1996}%
\special{pa 5549 2430}%
\special{fp}%
%
\special{pn 8}%
\special{pa 980 3050}%
\special{pa 1150 2150}%
\special{dt 0.045}%
\special{sh 1}%
\special{pa 1150 2150}%
\special{pa 1118 2212}%
\special{pa 1140 2202}%
\special{pa 1157 2219}%
\special{pa 1150 2150}%
\special{fp}%
\special{pa 1020 3030}%
\special{pa 1390 2420}%
\special{dt 0.045}%
\special{sh 1}%
\special{pa 1390 2420}%
\special{pa 1338 2467}%
\special{pa 1362 2466}%
\special{pa 1373 2487}%
\special{pa 1390 2420}%
\special{fp}%
%
\special{pn 8}%
\special{pa 2650 2820}%
\special{pa 2460 2310}%
\special{dt 0.045}%
\special{sh 1}%
\special{pa 2460 2310}%
\special{pa 2465 2379}%
\special{pa 2479 2360}%
\special{pa 2502 2365}%
\special{pa 2460 2310}%
\special{fp}%
\special{pa 2640 2800}%
\special{pa 2050 2340}%
\special{dt 0.045}%
\special{sh 1}%
\special{pa 2050 2340}%
\special{pa 2090 2397}%
\special{pa 2092 2373}%
\special{pa 2115 2365}%
\special{pa 2050 2340}%
\special{fp}%
\put(24.5000,-31.1000){\makebox(0,0)[lb]{$S^{p-1}\x D^{n+1-p}$}}%
\put(4.1000,-33.0000){\makebox(0,0)[lb]{$S^{n-p}\x D^p$}}%
\end{picture}%
\vskip2cm
Figure 2.

\np \input 3+.tex  
\vskip3cm
Figure 3.  $K_+$

\np \input 3-.tex  
\vskip3cm
Figure 3.  $K_-$

\np \input 30.tex  
\vskip3cm
Figure 3.  $K_0$

\np 
\unitlength 0.1in
\begin{picture}(56.10,43.70)(8.50,-44.30)
%
\special{pn 8}%
\special{ar 3510 320 560 250  0.0000000 6.2831853}%
%
\special{pn 20}%
\special{pa 3660 3950}%
\special{pa 3660 2730}%
\special{fp}%
%
\special{pn 20}%
\special{pa 3340 3970}%
\special{pa 3340 2740}%
\special{fp}%
%
\special{pn 20}%
\special{pa 3350 4000}%
\special{pa 3360 3970}%
\special{pa 3384 3949}%
\special{pa 3413 3935}%
\special{pa 3444 3926}%
\special{pa 3475 3921}%
\special{pa 3507 3920}%
\special{pa 3539 3923}%
\special{pa 3570 3930}%
\special{pa 3601 3940}%
\special{pa 3627 3958}%
\special{pa 3647 3983}%
\special{pa 3648 4014}%
\special{pa 3630 4040}%
\special{pa 3603 4058}%
\special{pa 3573 4070}%
\special{pa 3542 4077}%
\special{pa 3510 4080}%
\special{pa 3478 4079}%
\special{pa 3447 4075}%
\special{pa 3416 4066}%
\special{pa 3387 4053}%
\special{pa 3362 4032}%
\special{pa 3350 4003}%
\special{pa 3350 4000}%
\special{sp}%
%
\special{pn 20}%
\special{pa 3340 360}%
\special{pa 3340 1580}%
\special{fp}%
%
\special{pn 20}%
\special{pa 3660 340}%
\special{pa 3660 1570}%
\special{fp}%
%
\special{pn 20}%
\special{ar 3500 310 150 80  0.0000000 6.2831853}%
%
\special{pn 8}%
\special{ar 3510 4080 560 250  0.0000000 6.2831853}%
%
\special{pn 8}%
\special{pa 4080 330}%
\special{pa 4080 4080}%
\special{fp}%
%
\special{pn 8}%
\special{pa 2950 340}%
\special{pa 2950 4080}%
\special{fp}%
%
\special{pn 8}%
\special{ar 1410 330 560 250  0.0000000 6.2831853}%
%
\special{pn 8}%
\special{ar 1410 4090 560 250  0.0000000 6.2831853}%
%
\special{pn 8}%
\special{pa 1980 340}%
\special{pa 1980 4090}%
\special{fp}%
%
\special{pn 8}%
\special{pa 850 350}%
\special{pa 850 4090}%
\special{fp}%
%
\special{pn 8}%
\special{ar 5890 310 560 250  0.0000000 6.2831853}%
%
\special{pn 8}%
\special{ar 5890 4070 560 250  0.0000000 6.2831853}%
%
\special{pn 8}%
\special{pa 6460 320}%
\special{pa 6460 4070}%
\special{fp}%
%
\special{pn 8}%
\special{pa 5330 330}%
\special{pa 5330 4070}%
\special{fp}%
%
\special{pn 20}%
\special{ar 3500 1590 150 80  0.0000000 6.2831853}%
%
\special{pn 20}%
\special{ar 3500 2710 150 80  0.0000000 6.2831853}%
%
\special{pn 20}%
\special{ar 3520 2070 560 250  0.0000000 6.2831853}%
%
\special{pn 8}%
\special{pa 3660 1600}%
\special{pa 5730 1600}%
\special{dt 0.045}%
\special{pa 5730 1600}%
\special{pa 5729 1600}%
\special{dt 0.045}%
%
\special{pn 8}%
\special{pa 3680 2720}%
\special{pa 5750 2720}%
\special{dt 0.045}%
\special{pa 5750 2720}%
\special{pa 5749 2720}%
\special{dt 0.045}%
%
\special{pn 20}%
\special{ar 5850 1600 150 80  0.0000000 6.2831853}%
%
\special{pn 20}%
\special{ar 5850 2730 150 80  0.0000000 6.2831853}%
%
\special{pn 20}%
\special{pa 5690 1610}%
\special{pa 5690 2700}%
\special{fp}%
%
\special{pn 20}%
\special{pa 6020 1640}%
\special{pa 6020 2730}%
\special{fp}%
\put(12.0000,-46.0000){\makebox(0,0)[lb]{t=-0.5}}%
\put(32.0000,-46.0000){\makebox(0,0)[lb]{t=0}}%
\put(56.0000,-46.0000){\makebox(0,0)[lb]{t=0.5}}%
%
\special{pn 8}%
\special{ar 3510 320 560 250  0.0000000 6.2831853}%
%
\special{pn 8}%
\special{pa 2950 340}%
\special{pa 2950 4080}%
\special{fp}%
%
\special{pn 8}%
\special{pa 4080 330}%
\special{pa 4080 4080}%
\special{fp}%
%
\special{pn 8}%
\special{ar 3510 4080 560 250  0.0000000 6.2831853}%
%
\special{pn 20}%
\special{ar 3500 310 150 80  0.0000000 6.2831853}%
%
\special{pn 20}%
\special{pa 3660 340}%
\special{pa 3660 1570}%
\special{fp}%
%
\special{pn 20}%
\special{pa 3340 360}%
\special{pa 3340 1580}%
\special{fp}%
%
\special{pn 20}%
\special{pa 3350 4000}%
\special{pa 3360 3970}%
\special{pa 3384 3949}%
\special{pa 3413 3935}%
\special{pa 3444 3926}%
\special{pa 3475 3921}%
\special{pa 3507 3920}%
\special{pa 3539 3923}%
\special{pa 3570 3930}%
\special{pa 3601 3940}%
\special{pa 3627 3958}%
\special{pa 3647 3983}%
\special{pa 3648 4014}%
\special{pa 3630 4040}%
\special{pa 3603 4058}%
\special{pa 3573 4070}%
\special{pa 3542 4077}%
\special{pa 3510 4080}%
\special{pa 3478 4079}%
\special{pa 3447 4075}%
\special{pa 3416 4066}%
\special{pa 3387 4053}%
\special{pa 3362 4032}%
\special{pa 3350 4003}%
\special{pa 3350 4000}%
\special{sp}%
%
\special{pn 20}%
\special{pa 3340 3970}%
\special{pa 3340 2740}%
\special{fp}%
%
\special{pn 20}%
\special{pa 3660 3950}%
\special{pa 3660 2730}%
\special{fp}%
%
\special{pn 8}%
\special{ar 1410 330 560 250  0.0000000 6.2831853}%
%
\special{pn 8}%
\special{ar 1410 4090 560 250  0.0000000 6.2831853}%
%
\special{pn 8}%
\special{pa 1980 340}%
\special{pa 1980 4090}%
\special{fp}%
%
\special{pn 8}%
\special{pa 850 350}%
\special{pa 850 4090}%
\special{fp}%
%
\special{pn 8}%
\special{ar 5890 310 560 250  0.0000000 6.2831853}%
%
\special{pn 8}%
\special{ar 5890 4070 560 250  0.0000000 6.2831853}%
%
\special{pn 8}%
\special{pa 6460 320}%
\special{pa 6460 4070}%
\special{fp}%
%
\special{pn 8}%
\special{pa 5330 330}%
\special{pa 5330 4070}%
\special{fp}%
\end{picture}%
\vskip3cm
Figure 4.  $K_+$

\np 
\unitlength 0.1in
\begin{picture}(56.10,43.70)(8.50,-44.30)
%
\special{pn 8}%
\special{ar 3510 320 560 250  0.0000000 6.2831853}%
%
\special{pn 20}%
\special{pa 3660 3950}%
\special{pa 3660 2730}%
\special{fp}%
%
\special{pn 20}%
\special{pa 3340 3970}%
\special{pa 3340 2740}%
\special{fp}%
%
\special{pn 20}%
\special{pa 3350 4000}%
\special{pa 3360 3970}%
\special{pa 3384 3949}%
\special{pa 3413 3935}%
\special{pa 3444 3926}%
\special{pa 3475 3921}%
\special{pa 3507 3920}%
\special{pa 3539 3923}%
\special{pa 3570 3930}%
\special{pa 3601 3940}%
\special{pa 3627 3958}%
\special{pa 3647 3983}%
\special{pa 3648 4014}%
\special{pa 3630 4040}%
\special{pa 3603 4058}%
\special{pa 3573 4070}%
\special{pa 3542 4077}%
\special{pa 3510 4080}%
\special{pa 3478 4079}%
\special{pa 3447 4075}%
\special{pa 3416 4066}%
\special{pa 3387 4053}%
\special{pa 3362 4032}%
\special{pa 3350 4003}%
\special{pa 3350 4000}%
\special{sp}%
%
\special{pn 20}%
\special{pa 3340 360}%
\special{pa 3340 1580}%
\special{fp}%
%
\special{pn 20}%
\special{pa 3660 340}%
\special{pa 3660 1570}%
\special{fp}%
%
\special{pn 20}%
\special{ar 3500 310 150 80  0.0000000 6.2831853}%
%
\special{pn 8}%
\special{ar 3510 4080 560 250  0.0000000 6.2831853}%
%
\special{pn 8}%
\special{pa 4080 330}%
\special{pa 4080 4080}%
\special{fp}%
%
\special{pn 8}%
\special{pa 2950 340}%
\special{pa 2950 4080}%
\special{fp}%
%
\special{pn 8}%
\special{ar 1410 330 560 250  0.0000000 6.2831853}%
%
\special{pn 8}%
\special{ar 1410 4090 560 250  0.0000000 6.2831853}%
%
\special{pn 8}%
\special{pa 1980 340}%
\special{pa 1980 4090}%
\special{fp}%
%
\special{pn 8}%
\special{pa 850 350}%
\special{pa 850 4090}%
\special{fp}%
%
\special{pn 8}%
\special{ar 5890 310 560 250  0.0000000 6.2831853}%
%
\special{pn 8}%
\special{ar 5890 4070 560 250  0.0000000 6.2831853}%
%
\special{pn 8}%
\special{pa 6460 320}%
\special{pa 6460 4070}%
\special{fp}%
%
\special{pn 8}%
\special{pa 5330 330}%
\special{pa 5330 4070}%
\special{fp}%
%
\special{pn 20}%
\special{ar 3500 1590 150 80  0.0000000 6.2831853}%
%
\special{pn 20}%
\special{ar 3500 2710 150 80  0.0000000 6.2831853}%
%
\special{pn 20}%
\special{ar 3520 2070 560 250  0.0000000 6.2831853}%
\put(12.0000,-46.0000){\makebox(0,0)[lb]{t=-0.5}}%
\put(32.0000,-46.0000){\makebox(0,0)[lb]{t=0}}%
\put(56.0000,-46.0000){\makebox(0,0)[lb]{t=0.5}}%
%
\special{pn 8}%
\special{ar 3510 320 560 250  0.0000000 6.2831853}%
%
\special{pn 8}%
\special{pa 2950 340}%
\special{pa 2950 4080}%
\special{fp}%
%
\special{pn 8}%
\special{pa 4080 330}%
\special{pa 4080 4080}%
\special{fp}%
%
\special{pn 8}%
\special{ar 3510 4080 560 250  0.0000000 6.2831853}%
%
\special{pn 20}%
\special{ar 3500 310 150 80  0.0000000 6.2831853}%
%
\special{pn 20}%
\special{pa 3660 340}%
\special{pa 3660 1570}%
\special{fp}%
%
\special{pn 20}%
\special{pa 3340 360}%
\special{pa 3340 1580}%
\special{fp}%
%
\special{pn 20}%
\special{pa 3350 4000}%
\special{pa 3360 3970}%
\special{pa 3384 3949}%
\special{pa 3413 3935}%
\special{pa 3444 3926}%
\special{pa 3475 3921}%
\special{pa 3507 3920}%
\special{pa 3539 3923}%
\special{pa 3570 3930}%
\special{pa 3601 3940}%
\special{pa 3627 3958}%
\special{pa 3647 3983}%
\special{pa 3648 4014}%
\special{pa 3630 4040}%
\special{pa 3603 4058}%
\special{pa 3573 4070}%
\special{pa 3542 4077}%
\special{pa 3510 4080}%
\special{pa 3478 4079}%
\special{pa 3447 4075}%
\special{pa 3416 4066}%
\special{pa 3387 4053}%
\special{pa 3362 4032}%
\special{pa 3350 4003}%
\special{pa 3350 4000}%
\special{sp}%
%
\special{pn 20}%
\special{pa 3340 3970}%
\special{pa 3340 2740}%
\special{fp}%
%
\special{pn 20}%
\special{pa 3660 3950}%
\special{pa 3660 2730}%
\special{fp}%
%
\special{pn 8}%
\special{ar 1410 330 560 250  0.0000000 6.2831853}%
%
\special{pn 8}%
\special{ar 1410 4090 560 250  0.0000000 6.2831853}%
%
\special{pn 8}%
\special{pa 1980 340}%
\special{pa 1980 4090}%
\special{fp}%
%
\special{pn 8}%
\special{pa 850 350}%
\special{pa 850 4090}%
\special{fp}%
%
\special{pn 8}%
\special{ar 5890 310 560 250  0.0000000 6.2831853}%
%
\special{pn 8}%
\special{ar 5890 4070 560 250  0.0000000 6.2831853}%
%
\special{pn 8}%
\special{pa 6460 320}%
\special{pa 6460 4070}%
\special{fp}%
%
\special{pn 8}%
\special{pa 5330 330}%
\special{pa 5330 4070}%
\special{fp}%
%
\special{pn 8}%
\special{pa 3320 1590}%
\special{pa 1590 1590}%
\special{dt 0.045}%
\special{pa 1590 1590}%
\special{pa 1591 1590}%
\special{dt 0.045}%
%
\special{pn 8}%
\special{pa 3330 2740}%
\special{pa 1620 2740}%
\special{dt 0.045}%
\special{pa 1620 2740}%
\special{pa 1621 2740}%
\special{dt 0.045}%
%
\special{pn 20}%
\special{ar 1420 1590 150 80  0.0000000 6.2831853}%
%
\special{pn 20}%
\special{ar 1420 2720 150 80  0.0000000 6.2831853}%
%
\special{pn 20}%
\special{pa 1260 1600}%
\special{pa 1260 2690}%
\special{fp}%
%
\special{pn 20}%
\special{pa 1590 1630}%
\special{pa 1590 2720}%
\special{fp}%
\end{picture}%
\vskip3cm
Figure 4.  $K_-$

\np 
\unitlength 0.1in
\begin{picture}(56.10,47.00)(8.50,-47.60)
%
\special{pn 8}%
\special{ar 3510 320 560 250  0.0000000 6.2831853}%
%
\special{pn 20}%
\special{pa 3350 4000}%
\special{pa 3360 3970}%
\special{pa 3384 3949}%
\special{pa 3413 3935}%
\special{pa 3444 3926}%
\special{pa 3475 3921}%
\special{pa 3507 3920}%
\special{pa 3539 3923}%
\special{pa 3570 3930}%
\special{pa 3601 3940}%
\special{pa 3627 3958}%
\special{pa 3647 3983}%
\special{pa 3648 4014}%
\special{pa 3630 4040}%
\special{pa 3603 4058}%
\special{pa 3573 4070}%
\special{pa 3542 4077}%
\special{pa 3510 4080}%
\special{pa 3478 4079}%
\special{pa 3447 4075}%
\special{pa 3416 4066}%
\special{pa 3387 4053}%
\special{pa 3362 4032}%
\special{pa 3350 4003}%
\special{pa 3350 4000}%
\special{sp}%
%
\special{pn 20}%
\special{sh 0.600}%
\special{ar 3500 310 150 80  0.0000000 6.2831853}%
%
\special{pn 20}%
\special{ar 3510 4080 560 250  0.0000000 6.2831853}%
%
\special{pn 8}%
\special{pa 4080 330}%
\special{pa 4080 4080}%
\special{fp}%
%
\special{pn 8}%
\special{pa 2950 340}%
\special{pa 2950 4080}%
\special{fp}%
%
\special{pn 8}%
\special{ar 1410 330 560 250  0.0000000 6.2831853}%
%
\special{pn 8}%
\special{ar 1410 4090 560 250  0.0000000 6.2831853}%
%
\special{pn 8}%
\special{pa 1980 340}%
\special{pa 1980 4090}%
\special{fp}%
%
\special{pn 8}%
\special{pa 850 350}%
\special{pa 850 4090}%
\special{fp}%
%
\special{pn 8}%
\special{ar 5890 310 560 250  0.0000000 6.2831853}%
%
\special{pn 8}%
\special{ar 5890 4070 560 250  0.0000000 6.2831853}%
%
\special{pn 8}%
\special{pa 6460 320}%
\special{pa 6460 4070}%
\special{fp}%
%
\special{pn 8}%
\special{pa 5330 330}%
\special{pa 5330 4070}%
\special{fp}%
%
\special{pn 20}%
\special{ar 3520 2070 560 250  0.0000000 6.2831853}%
\put(12.0000,-46.0000){\makebox(0,0)[lb]{t=-0.5}}%
\put(32.0000,-46.0000){\makebox(0,0)[lb]{t=0}}%
\put(56.0000,-46.0000){\makebox(0,0)[lb]{t=0.5}}%
\put(30.3000,-49.3000){\makebox(0,0)[lb]{Figure 1.1}}%
%
\special{pn 8}%
\special{ar 3510 320 560 250  0.0000000 6.2831853}%
%
\special{pn 20}%
\special{pa 3350 4000}%
\special{pa 3360 3970}%
\special{pa 3384 3949}%
\special{pa 3413 3935}%
\special{pa 3444 3926}%
\special{pa 3475 3921}%
\special{pa 3507 3920}%
\special{pa 3539 3923}%
\special{pa 3570 3930}%
\special{pa 3601 3940}%
\special{pa 3627 3958}%
\special{pa 3647 3983}%
\special{pa 3648 4014}%
\special{pa 3630 4040}%
\special{pa 3603 4058}%
\special{pa 3573 4070}%
\special{pa 3542 4077}%
\special{pa 3510 4080}%
\special{pa 3478 4079}%
\special{pa 3447 4075}%
\special{pa 3416 4066}%
\special{pa 3387 4053}%
\special{pa 3362 4032}%
\special{pa 3350 4003}%
\special{pa 3350 4000}%
\special{sp}%
%
\special{pn 20}%
\special{sh 0.600}%
\special{ar 3500 310 150 80  0.0000000 6.2831853}%
%
\special{pn 8}%
\special{ar 3510 4080 560 250  0.0000000 6.2831853}%
%
\special{pn 8}%
\special{pa 4080 330}%
\special{pa 4080 4080}%
\special{fp}%
%
\special{pn 8}%
\special{pa 2950 340}%
\special{pa 2950 4080}%
\special{fp}%
%
\special{pn 8}%
\special{ar 1410 330 560 250  0.0000000 6.2831853}%
%
\special{pn 8}%
\special{ar 1410 4090 560 250  0.0000000 6.2831853}%
%
\special{pn 8}%
\special{pa 1980 340}%
\special{pa 1980 4090}%
\special{fp}%
%
\special{pn 8}%
\special{pa 850 350}%
\special{pa 850 4090}%
\special{fp}%
%
\special{pn 8}%
\special{ar 5890 310 560 250  0.0000000 6.2831853}%
%
\special{pn 8}%
\special{ar 5890 4070 560 250  0.0000000 6.2831853}%
%
\special{pn 8}%
\special{pa 6460 320}%
\special{pa 6460 4070}%
\special{fp}%
%
\special{pn 8}%
\special{pa 5330 330}%
\special{pa 5330 4070}%
\special{fp}%
%
\special{pn 20}%
\special{pa 3050 2220}%
\special{pa 3050 4200}%
\special{fp}%
%
\special{pn 20}%
\special{pa 3200 1910}%
\special{pa 3200 2100}%
\special{fp}%
%
\special{pn 20}%
\special{pa 3210 2740}%
\special{pa 3210 2580}%
\special{fp}%
%
\special{pn 20}%
\special{pa 3210 3090}%
\special{pa 3210 2920}%
\special{fp}%
%
\special{pn 20}%
\special{pa 3210 3420}%
\special{pa 3210 3320}%
\special{fp}%
%
\special{pn 20}%
\special{pa 3210 3880}%
\special{pa 3220 3590}%
\special{fp}%
%
\special{pn 20}%
\special{pa 3260 3890}%
\special{pa 3210 3880}%
\special{fp}%
%
\special{pn 20}%
\special{pa 3340 3950}%
\special{pa 3290 3930}%
\special{fp}%
%
\special{pn 20}%
\special{pa 3940 4130}%
\special{pa 3990 4180}%
\special{fp}%
%
\special{pn 20}%
\special{pa 3680 4000}%
\special{pa 3940 4130}%
\special{fp}%
%
\special{pn 20}%
\special{pa 4000 4190}%
\special{pa 4000 2170}%
\special{fp}%
%
\special{pn 20}%
\special{pa 3780 4250}%
\special{pa 3600 4070}%
\special{fp}%
%
\special{pn 20}%
\special{pa 3770 2300}%
\special{pa 3780 4250}%
\special{fp}%
%
\special{pn 20}%
\special{pa 3370 4290}%
\special{pa 3450 4090}%
\special{fp}%
%
\special{pn 20}%
\special{pa 3400 2320}%
\special{pa 3370 4290}%
\special{fp}%
%
\special{pn 20}%
\special{pa 3050 4180}%
\special{pa 3350 4060}%
\special{fp}%
%
\special{pn 20}%
\special{pa 3860 3920}%
\special{pa 3750 3930}%
\special{fp}%
%
\special{pn 20}%
\special{pa 3860 3650}%
\special{pa 3860 3920}%
\special{fp}%
%
\special{pn 20}%
\special{pa 3850 3320}%
\special{pa 3850 3530}%
\special{fp}%
%
\special{pn 20}%
\special{pa 3850 2770}%
\special{pa 3840 3060}%
\special{fp}%
%
\special{pn 20}%
\special{pa 3850 2370}%
\special{pa 3850 2630}%
\special{fp}%
%
\special{pn 20}%
\special{pa 3840 1890}%
\special{pa 3840 2080}%
\special{fp}%
%
\special{pn 20}%
\special{pa 3530 3850}%
\special{pa 3520 3900}%
\special{fp}%
%
\special{pn 20}%
\special{pa 3510 3530}%
\special{pa 3510 3820}%
\special{fp}%
%
\special{pn 20}%
\special{pa 3510 3140}%
\special{pa 3510 3310}%
\special{fp}%
%
\special{pn 20}%
\special{pa 3510 2780}%
\special{pa 3510 3000}%
\special{fp}%
%
\special{pn 20}%
\special{pa 3510 2390}%
\special{pa 3510 2640}%
\special{fp}%
%
\special{pn 20}%
\special{pa 3510 2180}%
\special{pa 3510 2240}%
\special{fp}%
%
\special{pn 20}%
\special{pa 3510 1860}%
\special{pa 3510 2040}%
\special{fp}%
\end{picture}%
\vskip3cm
Figure 4.  $K_0$

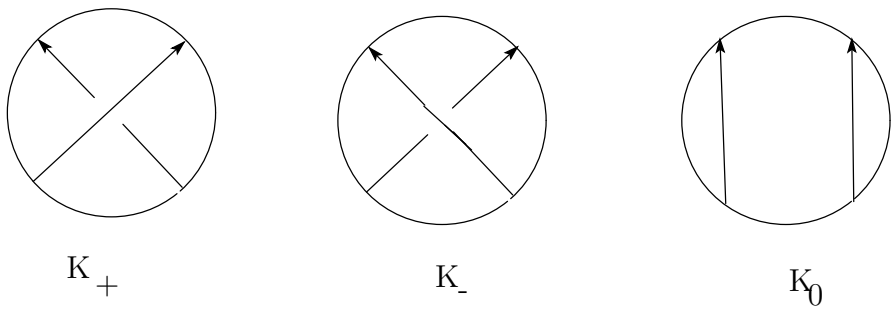
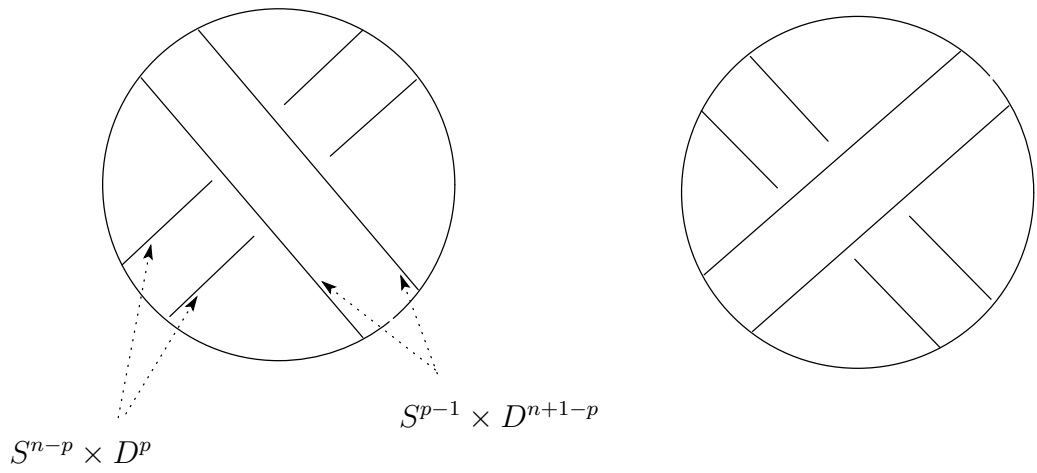
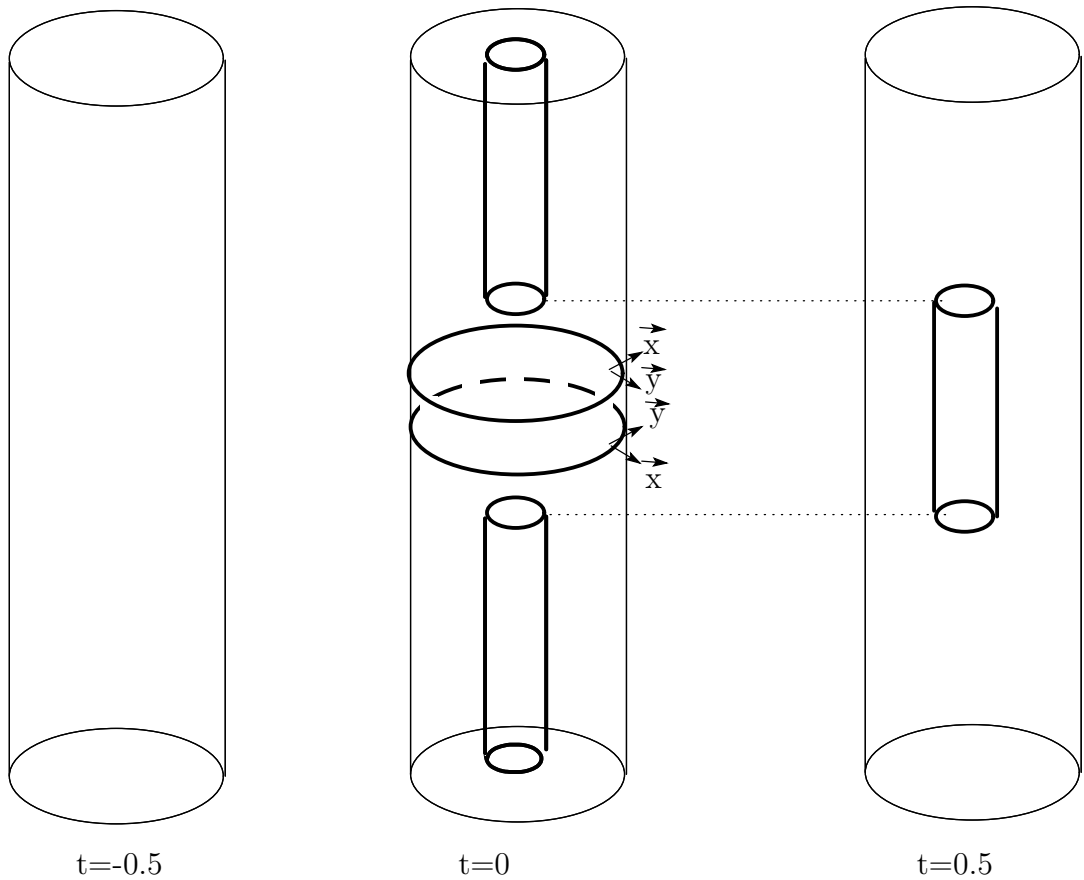
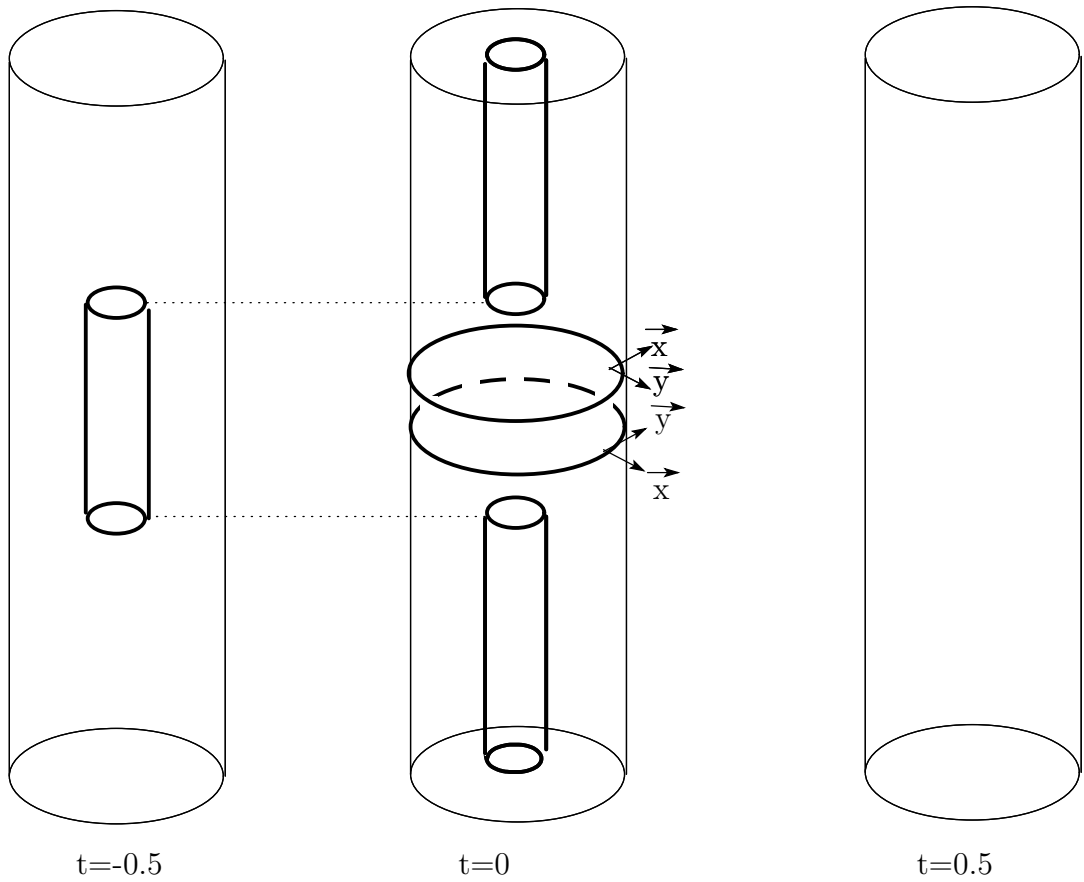
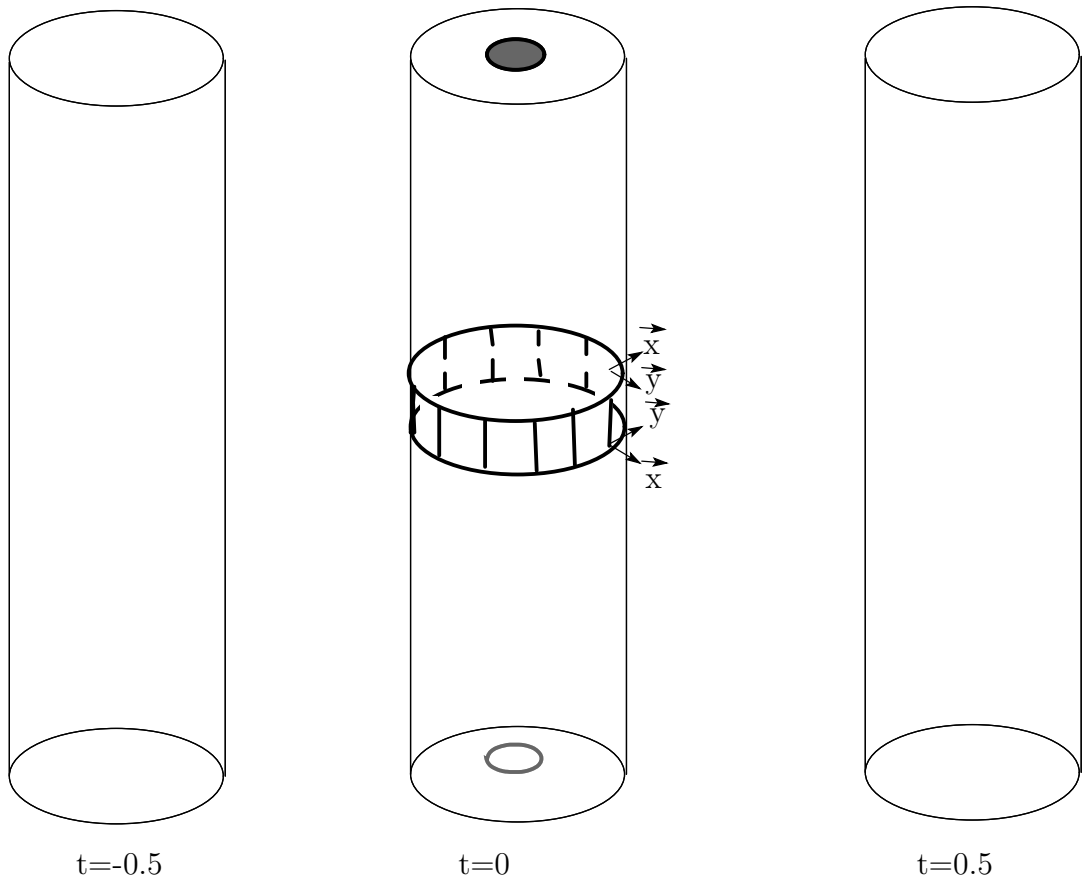
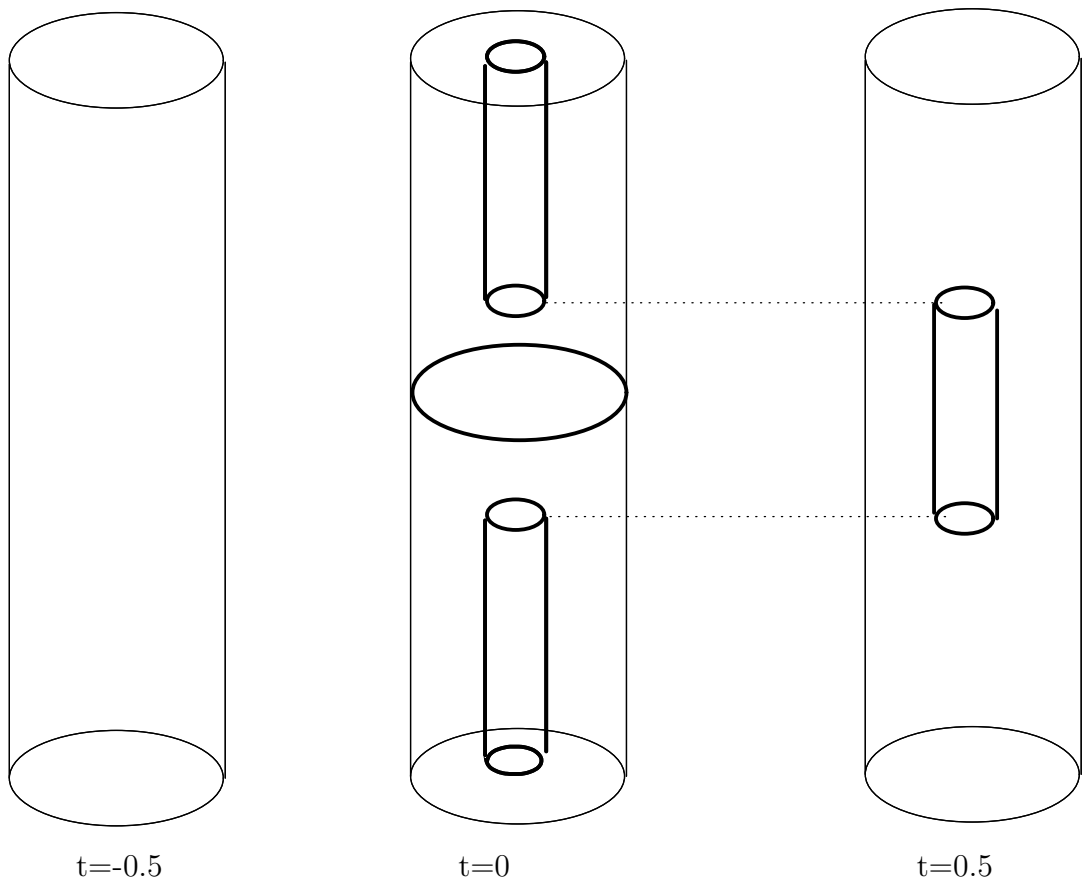
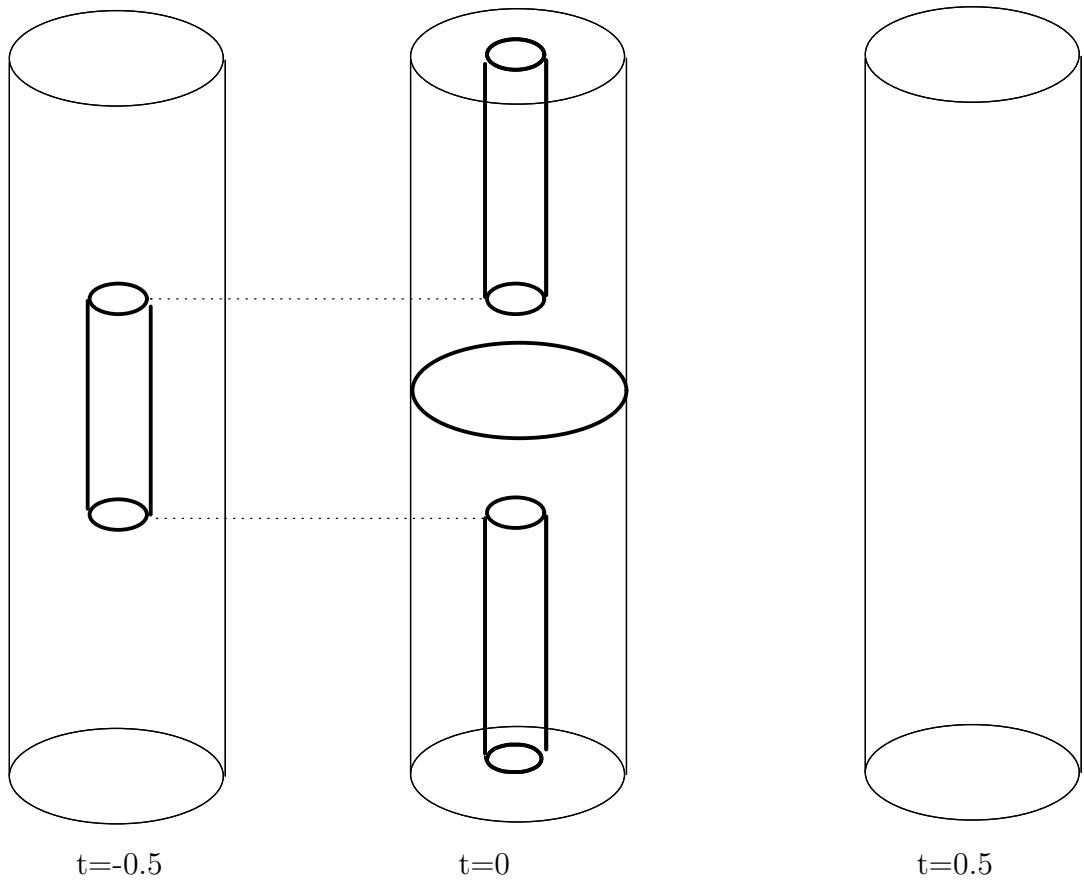
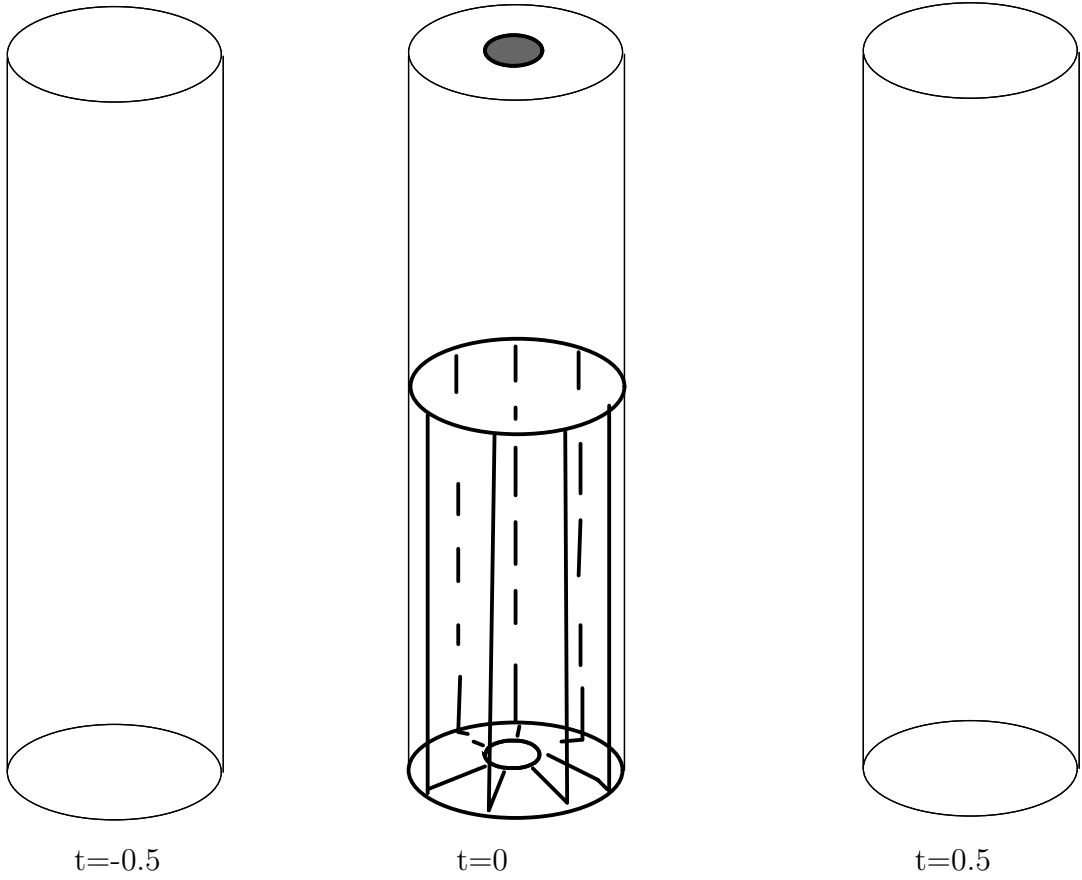
\end{document}